\documentclass[12pt]{amsart}
\pdfoutput=1
\usepackage{amsmath,amssymb,amscd,graphicx}
\usepackage[colorlinks=true, pdfstartview=FitV, linkcolor=blue, citecolor=blue, urlcolor=blue]{hyperref}
\usepackage{mathpazo}
\theoremstyle{plain}
\newtheorem*{nonumtheorem}{Theorem}

\newcommand{\tre}{\text{Re}}

\theoremstyle{definition}

\begin{document}
\title[Phase Statistics of the Riemann Zeros]{Notes on the Phase Statistics\\ of the Riemann Zeros}
\author{Jeffrey Stopple}
\address{Mathematics Department\\University of California, Santa Barbara\\Santa Barbara CA 93106}
\begin{abstract}
We numerically investigate, for zeros $\rho=1/2+i\gamma$, the statistics of the imaginary part of $\log(\zeta^\prime(1/2+i\gamma))$, computed by continuous variation along a vertical line from $\sigma=4$ to $4+i\gamma$ and then along a horizontal line to $1/2+i\gamma$.
\end{abstract}
\email{stopple@math.ucsb.edu}
%\keywords{Riemann hypothesis, de Bruijn-Newman constant, backward heat equation, Lehmer pair, random matrix theory}
%\subjclass[2010]{11M06, 11M50, 11Y35}

\maketitle

\section{Introduction}

\begin{figure}
\begin{center}
\includegraphics[scale=.75, viewport=0 0 450 450,clip]{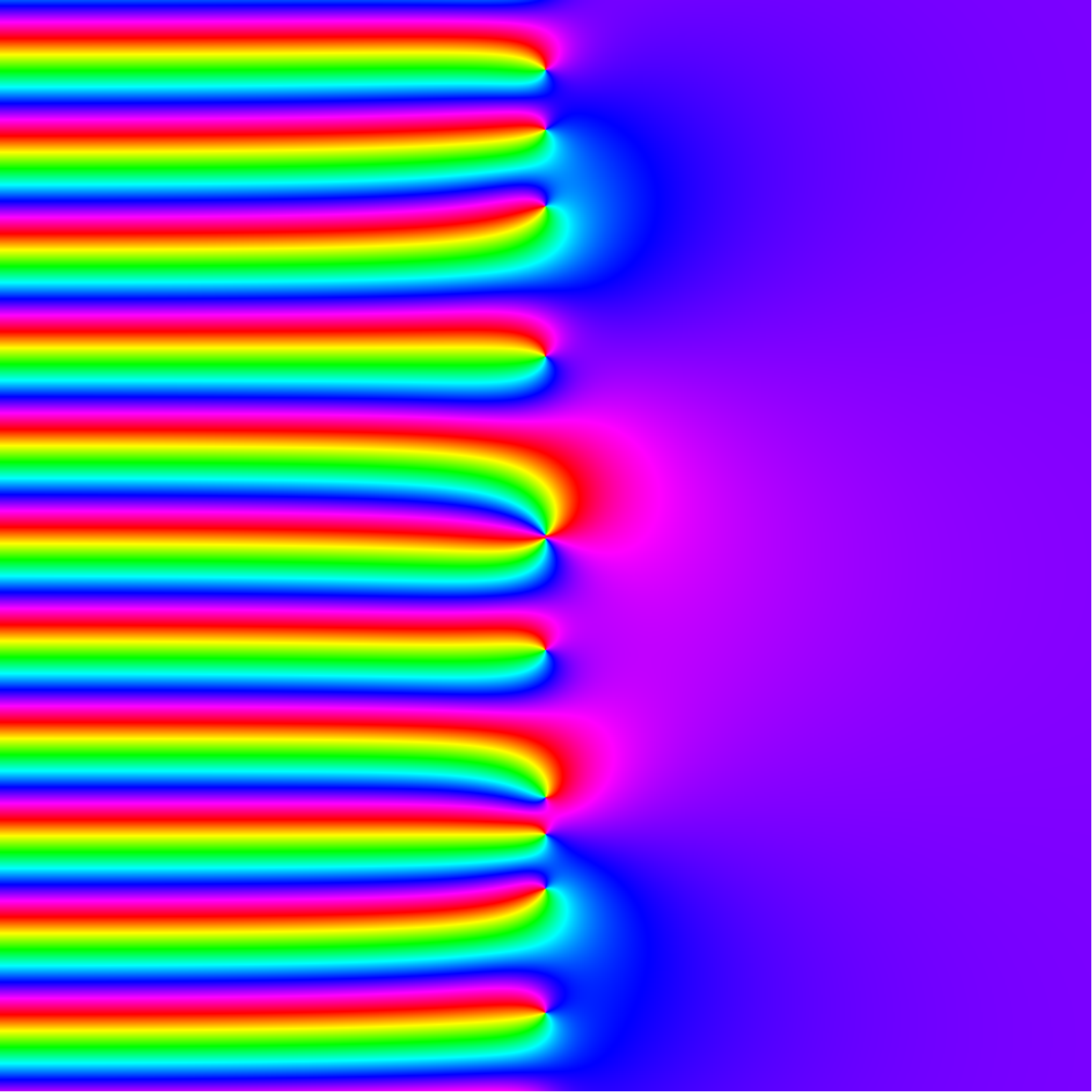}
\caption{$\arg(\zeta(s))$ for $-9/2\le \sigma\le11/2$, $7000\le t\le7010$.}\label{F:1}
\end{center}
\end{figure}

One popular way \cite{Wegert} of visualizing a complex function $w=f(z)$ is to plot $\arg(w)/2\pi$ interpreted as a color at each point $z$ in the domain.  This is easily implemented in  \emph{Mathematica}.  For the Riemann zeta function  the excitement is all near the critical strip: for $1\ll \tre(s)$, $\zeta(s)\approx 1$ and so the image is monochrome in that region.  For $\tre(s)\ll 0$, 
\[
\zeta(s)=\pi^{s-1/2}\frac{\Gamma((1-s)/2)}{\Gamma(s/2)}\zeta(1-s).
\]
Again $\zeta(1-s)\approx 1$. For bounded $\sigma$, as $t\to+\infty$,  Stirling's formula shows the argument of the remaining terms is asymptotic to $-t\log(t/2\pi)-t$, which means one sees very regular repeating horizontal bands of color.  Meanwhile, near a zero $\rho$ in the critical strip,
\[
\zeta(s)=\zeta^\prime(\rho)(s-\rho) +O(s-\rho)^2.
\]
Near $\rho$, the image corresponding to the function $s-\rho$ is just the color wheel with all the colors coming together at $s=\rho$.  Multiplying by $\zeta^\prime(\rho)$ locally scales the picture by $|\zeta^\prime(\rho)|$ and rotates it by $\arg(\zeta^\prime(\rho))$.  See Figure \ref{F:1} for an image with $7000\le t\le 7010$.  (What looks like a double zero is actually the first known example of a Lehmer pair near $t=7005.$)

Thus the argument of $\zeta^\prime(\rho)$ plays a significant role in the image, inspiring \href{https://mathoverflow.net/q/323852}{this} MathOverflow question.  In this paper we begin the numerical investigation of these by examining $5\cdot 10^6$ zeros with $2.63012\cdot 10^6\le \gamma\le 4.99238\cdot 10^6$.

\section{Selberg and Hejhal}

Starting first with $\zeta(s)$, an unpublished result of Selberg \cite[p.\ 310]{Titchmarsh} implies
\begin{nonumtheorem} Suitably normalized,
$\arg(\zeta(1/2+i t))$ converges in distribution over fixed ranges to a standard normal variable.  More precisely, for $\alpha<\beta$ we have
\begin{multline*}
\lim_{T\to\infty}\frac{1}{T}\mu(T\le t\le 2T\,|\,\alpha<\frac{\arg(\zeta(1/2+i t))}{\sqrt{\log\log(T)/2}}<\beta)=\\
\frac{1}{\sqrt{2\pi}}\int_\alpha^\beta \exp(-x^2/2)\, dx,
\end{multline*}
where $\mu$ is Lebesgue measure.
\end{nonumtheorem}

The argument is the imaginary part of $\log(\zeta(1/2+it))$, computed by continuous variation along a vertical line from, say, $\sigma=4$ to $4+it$ and then along a horizontal line to $1/2+it$.
Selberg's result actually covers the real part of the complex logarithm, $\log|\zeta(1/2+it)|$ as well.  

Returning to $\zeta^\prime$, for the real part of the logarithm we have the following generalization of Selberg's result, due to Hejhal \cite{Hejhal}:
\begin{nonumtheorem}
Assuming the Riemann Hypothesis and a technical condition on the spacing of zeros which is a weak consequence of the Montgomery Pair Correlation Conjecture,  then $\log |\zeta^\prime(1/2+i t)|$, suitably normalized,  converges in distribution over fixed ranges to a standard normal variable.  More precisely, for $\alpha<\beta$ we have

\begin{multline*}
 \lim_{N\to\infty}\left|\frac1N\left\{n:N\le n\le 2N,\alpha<\frac{\log\left|\frac{2\pi\zeta^\prime(1/2+i\gamma_n)}{\log(\gamma_n/2\pi)}\right|}{\sqrt{\log\log(N)/2}}<\beta\right \}\right|\\
 =
\frac{1}{\sqrt{2\pi}}\int_\alpha^\beta\exp(-x^2/2)\, dx
\end{multline*}
\end{nonumtheorem}
Here's my attempt to make an exposition of Hejhal's exposition \cite[p. 346]{Hejhal} of the basic idea behind the proof.  First some notation: With 
\[
\chi(s)=\pi^{s-1/2}\frac{\Gamma((1-s)/2)}{\Gamma(s/2)},
\]
define $\phi(s)$ by
\[
\chi(s)^{-1/2}=\exp(i\phi(s)), 
\]
so $\phi$ is real on the critical line.  Let $M$ be a large constant, $t$ an auxiliary random variable with $T\le t\le 2T$, and $W_t$ the 'window' 
\[
W_t=[t-M/\log T,t+M/\log T].
\]
Let $A(t)=t/2\pi\log\left(t/2\pi\right)-t/2\pi$. Let $x=A(u)$, and $\theta(u)=\phi(1/2+iu)$.  Let $P_t(x)$ be the polynomial approximation 
\[
P_t(x)=\Pi_{\gamma\in W_t}(x-A(\gamma)),
\]
and define $\Omega_t(u)$ to be the correction to the approximation so that
\[
\zeta(1/2+iu)=\exp(\Omega_t(u)-i\theta(u))P_t(x).
\]

Computing logarithmic derivative (in $u$, being careful with the chain rule) we see

$$
i\frac{\zeta^\prime(1/2+iu)}{\zeta(1/2+iu)}=\Omega_t^\prime(u)-i\theta^\prime(u)+\frac{P_t^\prime(x)}{P_t(x)}\cdot A^\prime(u)
$$

Rearranging gives

\[
i\frac{\zeta^\prime(1/2+iu)}{A^\prime(u)}=\zeta(1/2+iu)\cdot
\left(\frac{\Omega_t^\prime(u)}{A^\prime(u)}-i\frac{\theta^\prime(u)}{A^\prime(u)}+\frac{P_t^\prime(x)}{P_t(x)}\right).
\]

Hejhal  makes an estimate (see below) of the term in parenthesis on the right to argue that
\[
\frac{\log\left|\zeta(1/2+iu)\right|}{\sqrt{\log\log T}}\quad\text{ and }\quad\frac{\log\left|\zeta^\prime(1/2+iu)/A^\prime(u)\right|}{\sqrt{\log\log T}}
\]
are (in effect) the same random variable, and so Selberg's theorem applies.

For this estimate, Hejhal claims he and Bombieri showed previously that the total variation of $\Omega_t(u)$ on $W_t$ is $O_M(1)$ for 'most' $t$.  This means

$$
\int_{W_t}\left|\Omega^\prime_t(u)\right|du=O_M(1)
$$
for 'most' $t$, and so on 'most' windows $W_t$, $\left|\Omega^\prime_t(u)\right|=O_M(\log T)$.  

The above was the hard part; $\theta^\prime(u)/A^\prime(u)$ is elementary.  And $P_t^\prime(x)/P_t(x)=\sum_{\gamma\in W_t}1/(x-A(\gamma)),$ with average spacing between $A(\gamma)$ being 1 and the number of terms in the sum $O_M(1)$.  Hejhal argues heuristically that 

$$
\log\left|\frac{\Omega_t^\prime(u)}{A^\prime(u)}-i\frac{\theta^\prime(u)}{A^\prime(u)}+\frac{P_t^\prime(x)}{P_t(x)}\right|=O_M(1)
$$
except for a subset of small measure.  This completes Hejhal's estimate.

Could this heuristic be extended to the imaginary part of \newline  $\log\zeta^\prime(1/2+iu)$, defined (again) by continuous variation up the vertical line from $4$ to $4+iu$ and along the horizontal line from $4+iu$ to $1/2+iu$?  The challenge is that the estimates above depend on being inside a window of radius $M/\log T$.  We will see below that the variation along the vertical line is quite regular, so that presents no problem.  Along the horizontal line, between the real part $4$ and the real part $1/2+M/\log T$, the zeros of $\zeta^\prime$ should be infrequent and so the argument of $\zeta^\prime$ should be changing slowly for \lq most\rq values of $u$.

\section{Algorithm}
The first $10^7$ zeros of $\zeta(s)$ are implemented constants in \emph{Mathematica}, which can also, of course, easily compute derivatives numerically.  Evaluating the argument via continuous variation requires a little more effort.  For a zero $\rho=1/2+i\gamma$, the variation along the line from $4$  to $4+i\gamma$ is easy.  In this range
\[
\zeta^\prime(4+iy)=-\log(2)2^{-4-iy}-\sum_{n=3}^\infty\log(n)n^{-4-iy}.
\]
The first term has $|\log(2)2^{-4-iy}|=0.0433217$, while the tail is smaller, bounded by $0.025590$, and thus via Rouch\'e theorem $\zeta^\prime(4+iy)$ winds around the origin as many times as does $-\log(2)2^{-4-iy}$.  

Along the horizontal line $s=x+i\gamma$, $1/2<x\le4$, we only need to estimate $\zeta^\prime(s)$ very roughly, to determine when the argument increases by a multiple of $2\pi$.  Since we compute at many equally spaced points along the line, directly computing each derivative in \emph{Mathematica} is wasteful and slow.  Instead we make a table of values of $\zeta(s)$ at the equally spaced points, and use a variant of Richardson Interpolation \cite[5.7]{recipes} of the derivative:
\begin{multline*}
(-\zeta(s-3\Delta x)+9\zeta(s-2\Delta x)-45\zeta(s-\Delta x)+\\
45\zeta(s+\Delta x)-9\zeta(s+2\Delta x)+\zeta(s+3\Delta x))/60\Delta x=\\
\zeta^\prime(s)+O(\Delta x^6).
\end{multline*}
(In fact the next term of the $\Delta x$ series expansion of the left side  is\linebreak  $\zeta^{(7)}(s)\Delta x^6/140$.). The derivative $\zeta^\prime(1/2+i\gamma)$ is computed with the built-in \emph{Mathematica} implemetation.

\begin{table}
\begin{center}
\begin{tabular}{ c c } 
$k$ & $\Delta$  \\
\hline
$8546951$ & $0.00232317$  \\
$5042996$ & $0.00296997$  \\
$9857600$ & $0.00302828$  \\
$9675304$ & $0.00333645$  \\
$7279824$ & $0.00337943$  \\
$7498518$ & $0.00387655$  \\
$7060975$ & $0.00494052$  \\
$$ & $$ 
\end{tabular}
\caption{Small gaps between zeros $\rho_k$ of $\zeta(s)$, $5\cdot10^6\le k\le10^7$}\label{Ta:small}
\end{center}
\end{table}

The step size $\Delta x$ needs to give a sufficiently accurate result even when the horizontal line passes close to a zero of $\zeta^\prime(s)$.  Recalling that the zeros of $\zeta^\prime$ in $\tre(s)>1/2$ tend to be interspersed between the zeros of $\zeta(s)$ on the critical line, we looked for small gaps between successive zeros $\rho_k$, $\rho_{k+1}$ for $5\cdot 10^6\le k\le 10^7$.  With just seven exceptions (see Table \ref{Ta:small}), the gaps are all greater than $0.005$.  Based on this we choose for speed a step size $\Delta x$ of $0.0025$, accepting that a very small number of phases may be computed incorrectly.

\section{Data}
\begin{figure}
\begin{center}
\includegraphics[scale=1.25, viewport=0 0 425 160,clip]{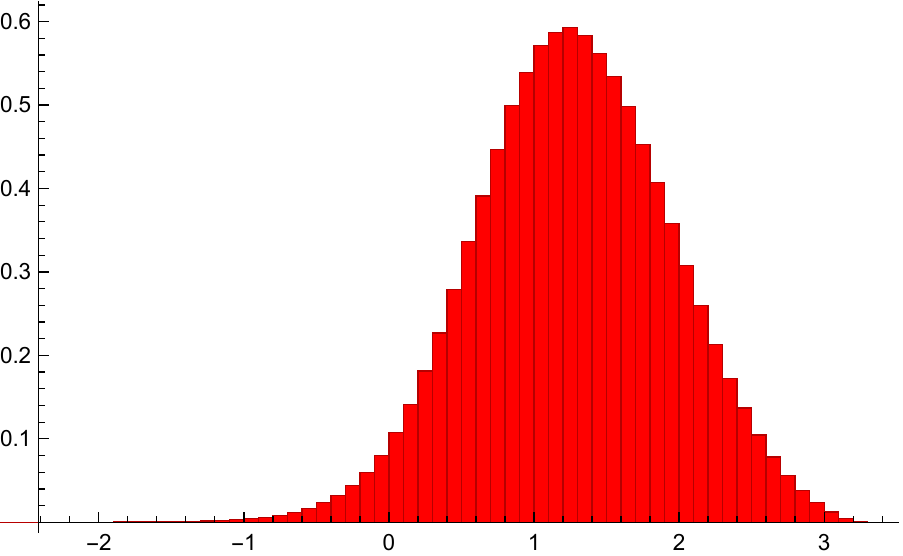}
\caption{$\log|\zeta^\prime(\rho)|$.}\label{F:3}
\end{center}
\end{figure}

Hiary and Odlyzko \cite{HO} have investigated Hejhal's theorem numerically, for data sets at much larger heights than we consider, and find the convergence rather slow.  They also observe a surplus of large values and a deficit of small values.    Since we have the data available, for completeness we include in Figure \ref{F:3} a histogram of values, for $5\cdot 10^6 \le k\le 10^7$, of 
\[
\frac{ \log|2\pi\zeta^\prime(\rho_k)/\log(\gamma_k/2\pi)|}{\sqrt{\log\log(5\cdot 10^6)}}.
\]

\begin{figure}
\begin{center}
\includegraphics[scale=1.25, viewport=0 0 425 160,clip]{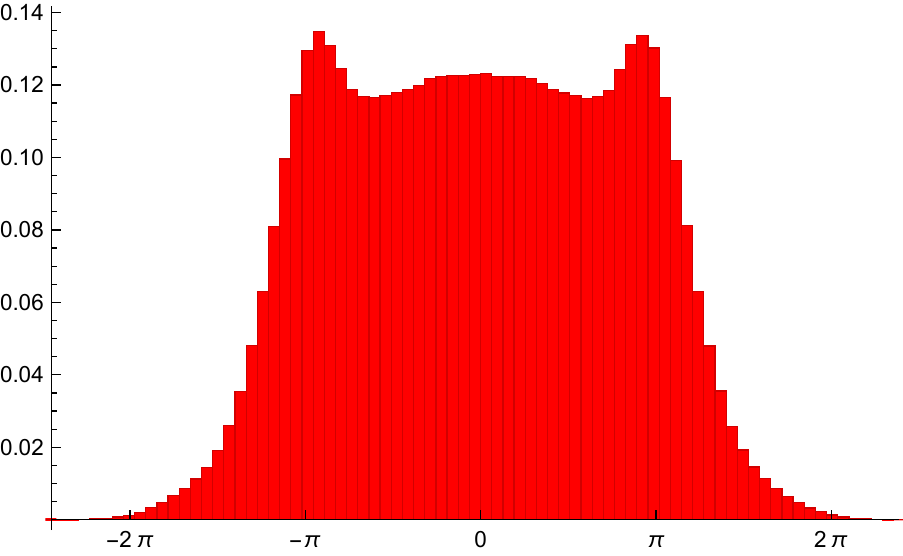}
\caption{$\arg(\zeta^\prime(\rho))$.}\label{F:2}
\end{center}
\end{figure}

Figure \ref{F:2} is the numerical investigation of the argument, the main goal of the paper.  For $5\cdot 10^6 \le k\le 10^7$, the histogram displays
\[
\frac{\arg\zeta^\prime(\rho_k)+\pi-\gamma_k\log 2}{\sqrt{\log\log(5\cdot 10^6)}} .
\]
\emph{Mathematica} computes the mean to be $-0.00043882$ and the standard deviation to be $2.47623$.  For what it is worth, the third through sixth moments were computed to be $0.00463054$, $76.8629$, $0.344781$, and $1333.96$ respectively.

\begin{figure}
\begin{center}
\includegraphics[scale=1.25, viewport=0 0 425 192,clip]{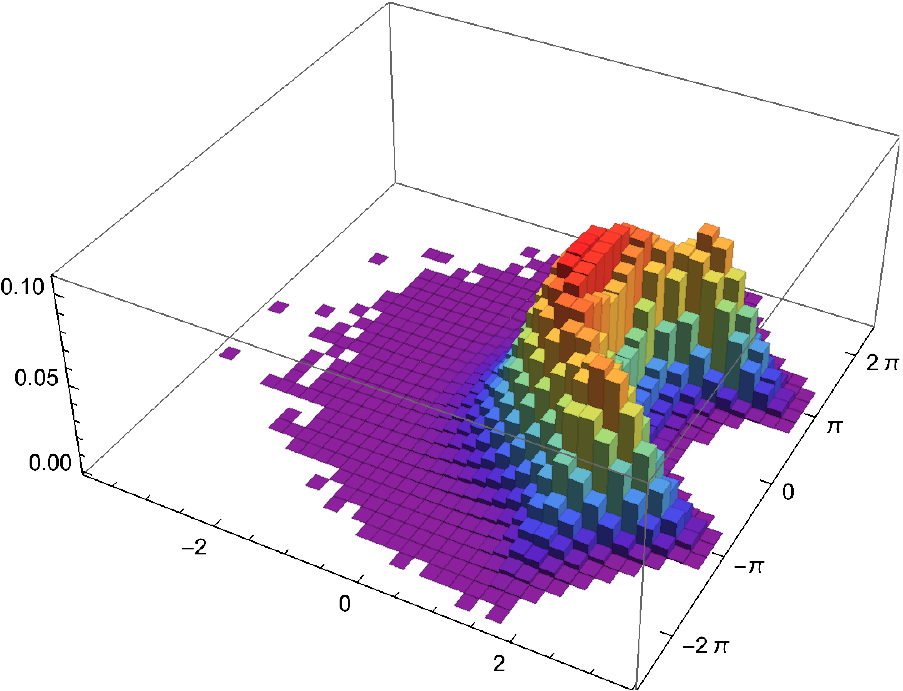}
\caption{$\log(\zeta^\prime(\rho))$.}\label{F:4}
\end{center}
\end{figure}

Figure \ref{F:4} shows both the real and imaginary parts of $\log \zeta^\prime(\rho_k)$.  Observe that the apparent surplus of examples with the imaginary part near $\pm\pi$ seems to correlate to the real part being positive and relatively large.

\section{Summary} Given the poor fit to a (mean $0$) Gaussian for the data in Figure \ref{F:3}, perhaps not much can be conjectured from the data in Figure \ref{F:2}, beyond that there \emph{is} a distribution for the argument computed by continuous variation.  In other words, the naive conjecture that the data are uniform in $(-\pi,\pi]$ appears to be incorrect.  We hope this paper inspires others with access to more computing power to investigate further.

\end{document}